\documentclass{article}
\usepackage{amsfonts}
\usepackage{amsmath}

\newtheorem{theorem}{Theorem}

\newtheorem{definition}[theorem]{Definition}

\newtheorem{proposition}[theorem]{Proposition}

\input{tcilatex}
\begin{document}

\title{Prehomogeneous geometries vs. \\
bundles with connections}
\author{Erc\"{u}ment H. Orta\c{c}gil}
\maketitle

\begin{abstract}
The below discussion is in three sections A, B, C, each section in two parts
I, II, I representing the standpoint of bundles with connections and II
representing the standpoint of prehomogeneous geometries (phg's). In A, our
object of study is a fibered manifold $\mathcal{E\rightarrow }M$, in B its $%
k $-jet bundle $J^{k}(\mathcal{E})\rightarrow M$ and in C a general
nonlinear PDE $\mathcal{H}^{k}\subset J^{k}(\mathcal{E}).$
\end{abstract}

\section{A.I.}

Let $\pi :\mathcal{E\rightarrow }M$ be a smooth fibered manifold fibering as 
$(x,y)\rightarrow x$ in coordinates ([P2], Chapter II). Suppose we want to
define a structure on $\pi :\mathcal{E\rightarrow }M$ which will enable us
to make the following decision: For any $q\in \mathcal{E}$ with $\pi (q)=p,$
there esists a \textit{unique} local section $s$ defined near $p$ satisfying 
$s(p)=q.$ This is the well known framework of connections: First, we
consider the $1$-jet bundle $J^{1}(\mathcal{E)\rightarrow E}$ whose fiber $%
J^{1}(\mathcal{E)}_{q}$ over $q\in \mathcal{E}$ is the set of all $1$-jets
of local sections $s$ of $\pi :\mathcal{E\rightarrow }M$ at $p=\pi (q).$ Now 
$J^{1}(\mathcal{E)\rightarrow E}$ is an affine bundle modeled on the vector
bundle $T^{\ast }(M)\otimes V(\mathcal{E})\rightarrow \mathcal{E}$ (where $%
T^{\ast }(M)\rightarrow M$ is pulled back to $\mathcal{E)}.$ In particular,
the fibers of $J^{1}(\mathcal{E)\rightarrow E}$ are contractible and $J^{1}(%
\mathcal{E)\rightarrow E}$ admits a global crossection $c:\mathcal{%
E\rightarrow }J^{1}(\mathcal{E)}.$

\begin{definition}
A crossection $c:\mathcal{E\rightarrow }J^{1}(\mathcal{E)}$ of the affine
bundle $J^{1}(\mathcal{E)\rightarrow E}$ is a connection on $\pi :\mathcal{%
E\rightarrow }M.$
\end{definition}

\textbf{Remark 1. }Some bundles $\pi :\mathcal{E\rightarrow }M$ may carry
some extra structure and it is natural to expect a connection to respect
this structure. With such a restriction a connection may not exist. We will
turn back to this issue below.

\bigskip

It turns out that the structure we search for is a flat connection, i.e., a
connection with vanishing curvature. Though well known, it is instructive to
see this construction explicitly in coordinates for our future purpose. So
let $q=(\overline{x},\overline{y})=(\overline{x}^{i},\overline{y}^{\alpha
})\in $ $U\times V,$ $1\leq i\leq n=\dim M,$ $1\leq \alpha \leq k=\dim 
\mathcal{E}_{p}$ and $s(x)$ a local section defined on $U$ satisfying $s(%
\overline{x})=\overline{y}.$ Now $(j^{1}s)(\overline{x})=(\overline{x}^{i},%
\overline{y}^{\alpha },\frac{\partial s^{\alpha }(\overline{x})}{\partial
x^{i}})\overset{def}{=}(\overline{x}^{i},\overline{y}^{\alpha },\overline{y}%
_{j}^{\alpha })$ and the values $\overline{y}_{i}^{\alpha }$ as $s(x)$
ranges over all such local section define the fiber $J^{1}(\mathcal{E)}_{q}.$
So a connection $c$ is of the form $(x^{i},y^{\alpha },c_{j}^{\alpha
}(x,y)). $ Now by the definition of $c,$ given $(\overline{x}^{i},\overline{y%
}^{\alpha },c_{j}^{\alpha }(\overline{x},\overline{y})),$ there exists a
local section $s(x)$ with $s(\overline{x})=\overline{y}$ satisfying

\begin{equation}
(\overline{x}^{i},\overline{y}^{\alpha },\frac{\partial s^{\alpha }(%
\overline{x})}{\partial x^{j}})=(\overline{x}^{i},\overline{y}^{\alpha
},c_{j}^{\alpha }(\overline{x},\overline{y}))
\end{equation}

Now $s(x)$ defines also the section $(x^{i},y^{\alpha },\frac{\partial
s^{\alpha }(x)}{\partial x^{j}}),$ $x\in Dom(s)$ but we need not have

\begin{equation}
(x^{i},y^{\alpha },\frac{\partial s^{\alpha }(x)}{\partial x^{j}}%
)=(x^{i},y^{\alpha },c_{j}^{\alpha }(x,y))\text{ \ \ \ }x\in Dom(s)
\end{equation}%
even though (2) holds at $x=\overline{x}$ by (1). The key fact here is that
the section in (1) depends on the point $(\overline{x},\overline{y})$ and
the same section may not work for all near points. We will write (2) shortly
as

\begin{equation}
\frac{\partial y^{\alpha }(x)}{\partial x^{j}}=c_{j}^{\alpha }(x,y)
\end{equation}

Now (3) is a first order system of PDE's with initial conditions, i.e.,
given the initial condition $(\overline{x},\overline{y}),$ we search for a
local section $y(x)$ solving (3) for $x\in U$ and satisfying $y(\overline{x}%
)=\overline{y}.$ This can be done if and only if the integrability
conditions of (3) are identically satisfied. To find these conditions, we
differentiate (3) with respect to $x^{k},$ substitute back from (3) and
alternate $j,k.$ The result is

\begin{equation}
\mathcal{R}_{jk}^{\alpha }(x,y)\overset{def}{=}\left[ \frac{\partial
c_{j}^{\alpha }(x,y)}{\partial x^{k}}+\frac{\partial c_{j}^{\alpha }(x,y)}{%
\partial y^{\beta }}c_{k}^{\beta }(x,y)\right] _{[jk]}=0
\end{equation}

It is easy to check that $\mathcal{R}$ is a $2$-form with values in the
vertical bundle $V(\mathcal{E})\rightarrow M,$ i.e., given $(\overline{x},%
\overline{y})\in \mathcal{E}$ and the two tangent vectors $\xi ,\eta $ at $%
\overline{x},$ $\mathcal{R}_{ab}^{\alpha }(\overline{x},\overline{y})\xi
^{a}\eta ^{b}$ is a vertical tangent vector at $(\overline{x},\overline{y}).$
We recall here the well known coordinatefree interpretation of $c$ and $%
\mathcal{R}:$ $c$ gives a decomposition of the tangent space $T_{q}(\mathcal{%
E})$ as a direct sum of horizontal and vertical components as

\begin{eqnarray}
T_{q}(\mathcal{E}) &=&H_{q}\oplus V(\mathcal{E})_{q} \\
(v^{i},v^{\alpha }) &=&(v^{i},v^{\alpha }-v^{s}c_{s}^{\alpha
})+(0,v^{s}c_{s}^{\alpha })  \notag
\end{eqnarray}%
and the horizontal distribution of dimension $n=\dim M$ defined by (5) is
integrable if and only if $\mathcal{R}=0.$

As the solution of our problem, we make the following

\begin{definition}
The unique local section $s(x)$ satisfying the initial condition $s(%
\overline{x})=\overline{y}$ is called a solution of the flat connection $c.$
\end{definition}

\section{A.II.}

Clearly $c(\mathcal{E})\subset J^{1}(\mathcal{E})$ by the definition of the
connection $c.$

\begin{definition}
The submanifold $c(\mathcal{E})\subset J^{1}(\mathcal{E})$ is a
prehomogeneous geometry (phg) of order one on the fibered manifold $\mathcal{%
E\rightarrow }M.$
\end{definition}

We observe that $c$ is viewed as part of the definition of a geometric
structure called a phg according to Definition 3. At first sight, a phg
seems to be nothing but a bundle with a connection on it but we will shortly
see that the situation is quite more subtle. In the language of [KS], [Gs],
[KLV], [P1], [P2], $c(\mathcal{E})$ is a nonlinear system of PDE's of order
one on $\mathcal{E\rightarrow }M.$ Note that the restriction of the jet
projection $\pi :J^{1}(\mathcal{E})\rightarrow \mathcal{E}$ gives a
bijection $\pi :c(\mathcal{E})\simeq \mathcal{E}$ so that $c(\mathcal{E}%
)\subset J^{1}(\mathcal{E})$ is a very special PDE. We recall ([P2], Chapter
3) that the prolongation $\varrho (c(\mathcal{E}))$ of $c(\mathcal{E})$ is
defined by

\begin{equation}
\varrho (c(\mathcal{E}))\overset{def}{=}J^{1}(c(\mathcal{E}))\cap J^{2}(%
\mathcal{E})\subset J^{1}J^{1}(\mathcal{E})
\end{equation}

Since $J^{1}(c(\mathcal{E}))$ surjects onto $c(\mathcal{E})$ by the jet
projection, (6) gives the map

\begin{equation}
\varrho (c(\mathcal{E}))\rightarrow c(\mathcal{E})
\end{equation}%
which need not be onto.

\begin{proposition}
\textit{(7) is onto if and only if }$\mathcal{R}=0.$
\end{proposition}

To see this, we first recall the canonical injection ([P2], Chapter 2, Lemma
4))

\begin{equation}
J^{2}(\mathcal{E})\subset J^{1}J^{1}(\mathcal{E})
\end{equation}%
Now $(x^{i},y^{\alpha },y_{j}^{\alpha },y_{jk}^{\alpha })$ are coordinates
on $J^{2}(\mathcal{E}),$ $(x^{i},y^{\alpha },y_{j}^{\alpha }\mid
y_{,k}^{\alpha },y_{j,k}^{\alpha })$ are coordinates on $J^{1}J^{1}(\mathcal{%
E})$ and (8) is given by

\begin{equation}
(x^{i},y^{\alpha },y_{j}^{\alpha },y_{jk}^{\alpha })\rightarrow
(x^{i},y^{\alpha },y_{i}^{\alpha }\mid y_{j}^{\alpha },y_{jk}^{\alpha })
\end{equation}%
i.e., the image of $J^{2}(\mathcal{E})$ in $J^{1}J^{1}(\mathcal{E})$ is
defined by the equations $y_{,k}^{\alpha }=y_{k}^{\alpha },$ $%
y_{j,k}^{\alpha }=y_{jk}^{\alpha }.$ Clearly $J^{1}(c(\mathcal{E}))\subset
J^{1}J^{1}(\mathcal{E})$ since $c(\mathcal{E})\subset J^{1}(\mathcal{E})$
and therefore $J^{1}(c(\mathcal{E}))$ is given by

\begin{equation}
J^{1}(c(\mathcal{E})):(x^{i},y^{\alpha },c_{j}^{\alpha }(x,y)\mid
y_{,k}^{\alpha },\frac{\partial c_{j}^{\alpha }(x,y)}{\partial x^{k}}+\frac{%
\partial c_{j}^{\alpha }(x,y)}{\partial y^{\beta }}y_{,k}^{\beta })
\end{equation}%
where $y_{,k}^{\alpha }=\frac{\partial y^{\alpha }}{\partial x^{k}}.$ It
follows from (10) that $\varrho (c(\mathcal{E}))=$ $J^{1}(c(\mathcal{E}%
))\cap J^{2}(\mathcal{E})$ surjects onto $c(\mathcal{E})$ if and only if the
substitution $y_{,k}^{\alpha }=c_{j}^{\alpha }(x,y)$ makes $y_{jk}^{\alpha }%
\overset{def}{=}$\ $\frac{\partial c_{j}^{\alpha }(x,y)}{\partial x^{k}}+%
\frac{\partial c_{j}^{\alpha }(x,y)}{\partial y^{\beta }}y_{,k}^{\beta }$\
symmetric in $j,k,$\ i.e, if and only if $\mathcal{R}=0.$

\bigskip

\textbf{Remark 2. }A very important fiber bundle in geometry is the
following: Consider the projection on the first factor $\pi :$ $M\times
M\rightarrow M$ whose local sections are local maps $M\rightarrow M.$ We
restrict our attention to sections which are local diffeomorphisms. With
this agreed, $J^{1}(M\times M)$ is $1$-jets of local diffeomorphisms called $%
1$-arrows in [Or1]. Now a connection $c$ (denoted by $\varepsilon $ in
[Or1]) assigns to the pair $(x,y)\in M\times M$ a $1$-arrow from $x$ to $y.$
These $1$-arrows have the structure of a groupoid and it is natural to
assume that connections respect this groupoid structure. Now such
connections exist if and only if $M$ is parallelizable and the first two
parts of [Or1] is devoted to the study of these phg's of order one (called
order zero in [Or1]. It turns out that the theory of Lie groups and Lie
algebras can be deduced from the study of these particular phg's. The name
prehomogeneous is motivated by this case because these structures make $M$
locally homogeneous (or endow $M$ with a geometric structure in the sense of
[T], [Gm]) when their curvatures vanish.

\section{B.I.}

Now we specialize our fiber bundles. Let $\pi :\mathcal{E}\rightarrow M$ be
a fibered manifold and our object of study in this section is the $k$-jet
bundle $J^{k}(\mathcal{E})\rightarrow M$ whose fiber over $q\in \mathcal{E}$
is $k$-jets of local sections $s(x)$ of $\mathcal{E}\rightarrow M$ with $%
s(p)=q.$ Now $J^{1}(J^{k}(\mathcal{E}))\rightarrow $ $J^{k}(\mathcal{E})$ is
an affine bundle modelled on $T^{\ast }(M)\otimes J^{k}(V(\mathcal{E}%
))\simeq T^{\ast }(M)\otimes V(J^{k}(\mathcal{E}))$ ([P2], Propositions 3,
5) and therefore admits a global crossection $c.$ If the curvature of this
connection vanishes, then for any $r\in J^{k}(\mathcal{E}),$ there exists a
unique local section $s_{k}(x)$ of $J^{k}(\mathcal{E})\rightarrow M$ passing
through $r.$ It is very crucial to observe that this section is not
necessarily holonomic, i.e., it may not satisfy $j^{k}(s)(p)=r$ where $p$ is
the projection of $r.$ However, we are interested in holonomic sections
because it is such sections which are true solutions of PDE's. To see this
point clearly, it is instructive to look at $c$ in coordinates. For
simplicity of notation, we assume $k=1.$ Now we have

\begin{eqnarray}
J^{1}(\mathcal{E}) &:&(x^{i},y^{\alpha },y_{j}^{\alpha })  \notag \\
J^{1}J^{1}(\mathcal{E)} &:&(x^{i},y^{\alpha },y_{j}^{\alpha }\mid
y_{,j}^{\alpha },y_{k,j}^{\alpha }) \\
c &:&(x^{i},y^{\alpha },y_{j}^{\alpha }\mid c_{,j}^{\alpha
}(x,y_{1}),c_{k,j}^{\alpha }(x,y_{1})  \notag
\end{eqnarray}%
where $(x,y_{1})\overset{def}{=}(x^{s},y^{\beta },y_{m}^{\beta }).$ If $r=(%
\overline{x}^{i},\overline{y}^{\alpha },\overline{y}_{j}^{\alpha })\in J^{1}(%
\mathcal{E}),$ a section $s$ of $J^{1}(\mathcal{E})\rightarrow M$ passing
through $r$ is of the form $s(x)=(x^{i},y^{\alpha }(x),y_{j}^{\alpha }(x))$
where $y^{\alpha }(\overline{x})=\overline{y}^{\alpha },y_{j}^{\alpha }(%
\overline{x})=\overline{y}_{j}^{\alpha }$ and the PDE for such sections is

\begin{eqnarray}
\frac{\partial y^{\alpha }}{\partial x^{j}} &=&c_{,j}^{\alpha }(x,y_{1}) \\
\frac{\partial y_{j}^{\alpha }(x)}{\partial x^{k}} &=&c_{k,j}^{\alpha
}(x,y_{1})  \notag
\end{eqnarray}

The integrability conditions of (12) are abtained by differentiating (12),
substituting back from (12) and alternating. They are given by%
\begin{equation*}
\mathcal{R}_{rj}^{\alpha }(x,y_{1})\overset{def}{=}
\end{equation*}

\begin{equation}
\left[ \frac{\partial c_{,j}^{\alpha }(x,y_{1})}{\partial x^{r}}+\frac{%
\partial c_{,j}^{\alpha }(x,y_{1})}{\partial y^{\beta }}c_{,r}^{\beta
}(x,y_{1})+\frac{\partial c_{,j}^{\alpha }(x,y_{1})}{\partial y_{a}^{\beta }}%
c_{r,a}^{\beta }(x,y_{1})\right] _{[rj]}=0  \notag
\end{equation}

\begin{equation}
\mathcal{R}_{rj,k}^{\alpha }(x,y_{1})\overset{def}{=}
\end{equation}%
\qquad

\begin{equation*}
\left[ \frac{\partial c_{k,j}^{\alpha }(x,y_{1})}{\partial x^{r}}+\frac{%
\partial c_{k,j}^{\alpha }(x,y_{1})}{\partial y^{\beta }}c_{,r}^{\beta
}(x,y_{1})+\frac{\partial c_{k,j}^{\alpha }(x,y_{1})}{\partial y_{a}^{\beta }%
}c_{r,a}^{\beta }(x,y_{1})\right] _{[rj]}=0
\end{equation*}%
and as in (4), the total curvature $\mathcal{R}(x,y_{1})\overset{def}{=}$ $(%
\mathcal{R}_{rj}^{\alpha }(x,y_{1}),\mathcal{R}_{rj,k}^{\alpha }(x,y_{1}))$
defined on $J^{1}(\mathcal{E})$ is easily seen to be a $2$-form on $M$ with
values in the vertical bundle $V(J^{1}(\mathcal{E}))\simeq J^{1}(V(\mathcal{E%
})).$ The computation for general $J^{k}(\mathcal{E})$ is similar and the
formulas can be written in compact form using the multi index notation $%
(x^{i},y^{\alpha },y_{j_{1}}^{\alpha },y_{j_{1}j_{2}}^{\alpha
},...,y_{j_{1}...j_{k}}^{\alpha })=(x^{i},y_{\mu }^{\alpha }),$ $0\leq
\left\vert \mu \right\vert \leq k,$ $=(x,y_{k}).$ Now if $\mathcal{R=}0,$
then $s(x)=(x^{i},y^{\alpha }(x),y_{j}^{\alpha }(x))$ satisfying (12) and
the above initial condition is unique but this solution does not necessarily
satisfy $y_{j}^{\alpha }(x)=\frac{\partial y^{\alpha }(x)}{\partial x^{j}},$
i.e., it is not necessarily holonomic. Thus the following question is of
fundamental importance.

\bigskip

\textbf{Q. }What property must the flat connection posess so that its
solutions are holonomic?

\section{B.II.}

We recall that $J^{k+1}(\mathcal{E})\rightarrow J^{k}(\mathcal{E})$ is an
affine bundle modelled on $S^{k+1}(T^{\ast }(M))\otimes V(\mathcal{E})$ and
therefore admits a global section $c.$ Now (8) generalizes to the canonical
inclusion $J^{k+1}(\mathcal{E})\subset J^{1}(J^{k}(\mathcal{E}))$ ([P2],
Chapter 2, Lemma 4). Therefore $c$ is also a section of $J^{1}(J^{k}(%
\mathcal{E}))\rightarrow J^{k}(\mathcal{E}),$ i.e., it is a connection on $%
J^{k}(\mathcal{E}).$

\begin{definition}
Connections on $J^{k}(\mathcal{E})\rightarrow M,$ i.e., crossections of $%
J^{1}(J^{k}(\mathcal{E}))\rightarrow J^{k}(\mathcal{E}),$ which arise from
the canonical inclusion $J^{k+1}(\mathcal{E})\subset J^{1}(J^{k}(\mathcal{E}%
))$ are geometric.
\end{definition}

\textbf{Remark 3. }All connections in A.I. are clearly geometric. For the
groupoids $J^{k}(M\times M),$ geometric connections exist only for $k=0$
(absolute parallelizm), $k=2,$ $\dim M$ arbitrary (affine) and $k=3,$ $\dim
M=1$ (projective) because jet groups split only for these values.

\bigskip

The following proposition answers \textbf{Q.}

\begin{proposition}
Solutions of flat geometric connections are holonomic.
\end{proposition}

To see what is involved in Proposition 6, let us look at the case $k=1$
above. It follows from (9) that a geometric connection $c$ satisfies

\begin{eqnarray}
c_{,j}^{\alpha }(x,y_{1}) &=&y_{j}^{\alpha } \\
c_{k,j}^{\alpha }(x,y_{1}) &=&y_{kj}^{\alpha }  \notag
\end{eqnarray}

In particular note that $c_{k,j}^{\alpha }(x,y_{1})$ is symmetric in $k,j.$

\bigskip

\textbf{Remark 4. }The above symmetry condition is interpreted as
torsionfreeness of $c$ which is meaningful only for $k=1$ ([Or3]). It is not
easy to detect this fact for large $k$ using the sophisticated formulation
(5).

\bigskip

Therefore $s(x)=(x^{i},y^{\alpha }(x),y_{j}^{\alpha }(x))$ solves

\begin{eqnarray}
\frac{\partial y^{\alpha }}{\partial x^{j}} &=&y_{j}^{\alpha }(x) \\
\frac{\partial y_{j}^{\alpha }(x)}{\partial x^{k}} &=&c_{kj}^{\alpha
}(x,y_{1})  \notag
\end{eqnarray}

The integrability conditions of (15) are given \textit{only }by the second
formula of (13) and the first equation of (15) shows that the solution is
holonomic, in fact, $\frac{\partial ^{2}y^{\alpha }}{\partial x^{k}\partial
x^{j}}=c_{kj}^{\alpha }.$ For general $k,$ all lower order integrability
conditions in (15) are satisfied except the top order which gives the
curvature. This is due to the fact that we work with $J^{k}(\mathcal{E})$
and not with a general nonlinear PDE $\mathcal{H}^{k}\subset J^{k}(\mathcal{E%
})$ which we will come to below.

\begin{definition}
Let $c$ be a geometric connection on $J^{k}(\mathcal{E}).$ Then the
nonlinear PDE $c(J^{k}(\mathcal{E}))\subset J^{k+1}(\mathcal{E})$ is a phg
of order $k+1.$
\end{definition}

As before, note that $c$ is now part of the definition of a geometric
structure. The following proposition, whose proof is identical to the
computation at the end of A.II, should not come as a surprise.

\begin{proposition}
The following are equivalent for a geometric connection $c.$

i) $\mathcal{R}=0$

ii) The prolongation map $\varrho (c(J^{k}(\mathcal{E})))\rightarrow $ $%
c(J^{k}(\mathcal{E}))$ is surjective.
\end{proposition}

As long as we work with connections on $J^{k}(\mathcal{E}),$ geometric ones
are clearly preferable in view of Proposition 6. From this standpoint, a phg
emerges as an alternative to a bundle with connection and it seems at this
stage that we are at a point of departure from connections. Furthermore,
Proposition 8, ii) hints at the possibility of defining the curvature of a
phg independently of the curvature of any connection. Below we will see that
this departure is not only a technical possibility, but also a very natural
alternative choice.

\bigskip

Below we will specialize to a fibered submanifold $\mathcal{H}^{k}\subset
J^{k}(\mathcal{E}),$ i.e., a general nonlinear system of PDE's. This is a
vast class and a highly sophisticated machinery is developed in [S], [KS],
[Gm], [KLV], [Ol1], [P1], [P2] (see also the references in these works) to
study the formal properties of these PDE's. Unfortunately this machinery is
not widely known and used among the differential geometers. The reason for
this state of affairs, we believe, is that this theory has not been able to
produce significant global examples and their invariants for large $k$ and
the well known examples for small $k$ (like affine, riemannian, projective,
conformal structures) can be thoroughly studied by the classical methods
without this machinery. In [Or4], we gave a completely elementary method to
construct phg's (see Definition 9 below) for arbitrarily large $k$ using the
irreducable representations of semi simple Lie algebras. For instance, using 
$SL(2,\mathbb{R}),$ it is possible to construct phg's for arbitrarily large $%
k$ on all $4$-manifolds. The question whether there exists a flat one among
all these phg's is quite relevant as Poincare Conjecture is a problem of
this type (see Remark 1 and [Or1], Chapter 13) Motivated by these facts, we
believe that jet theory has the potential to open new horizons in geometry
for large $k$ and has something new to offer also in the study of the well
known classical structures.

\section{C.I.}

The tale of connections is the same: We choose a crossection of $J^{1}(%
\mathcal{H}^{k})\rightarrow \mathcal{H}^{k},$ define its curvature which
vanishes if and only if there exist unique not necessarily holonomic
sections of $\mathcal{H}^{k}\rightarrow M$ satisfying the given initial
conditions.

\section{C.II}

The restriction of $J^{k+1}(\mathcal{E})\rightarrow J^{k}(\mathcal{E})$ to $%
\mathcal{H}^{k}\subset J^{k}(\mathcal{E})$ gives the fibered manifold $%
J^{k+1}(\mathcal{E})_{\mid \mathcal{H}^{k}}\rightarrow \mathcal{H}^{k}$ with
contractible fibers and therefore admits a crossection $\varepsilon .$

\begin{definition}
The nonlinear PDE $\mathcal{H}^{k}\simeq \varepsilon (\mathcal{H}%
^{k})\subset J^{k+1}(\mathcal{E})$ of order $k+1$ is a phg of order $k+1.$
\end{definition}

Therefore any nonlinear PDE $\mathcal{H}^{k}\subset J^{k}(\mathcal{E})$ of
order $k$ can be extended to a phg $\varepsilon (\mathcal{H}^{k})\subset
J^{k+1}(\mathcal{E})$ of order $k+1.$ \textit{The arbitrariness in the
choice of }$\varepsilon $\textit{\ is quite subtle and reflects the
"modelfreeness of the geometry" ([B], [Or1]).} Now even though $J^{k+1}(%
\mathcal{E})\subset J^{1}(J^{k}(\mathcal{E}))$ we do not necessarily have $%
\varepsilon (\mathcal{H}^{k})\subset J^{1}(\mathcal{H}^{k}),$ i.e., $%
\varepsilon $\textit{\ need not be a connection !! }Now $\varrho (\mathcal{H}%
^{k})=J^{1}(\mathcal{H}^{k})\cap J^{k+1}(\mathcal{E)}$ (we should not
confuse $\varrho (\mathcal{H}^{k})$ and $\varrho (\varepsilon (\mathcal{H}%
^{k}))=J^{1}(\varepsilon (\mathcal{H}^{k}))\cap J^{k+2}(\mathcal{E}))).$
Therefore $\varepsilon (\mathcal{H}^{k})\subset J^{1}(\mathcal{H}^{k}),$
i.e., $\varepsilon $ is a connection $\Leftrightarrow $ $\varepsilon (%
\mathcal{H}^{k})\subset \varrho (\mathcal{H}^{k}).$ If we stay on the side
of I, the question is whether we can choose $\varepsilon $ to be a
connection, and if so, whether it is unique,...etc (the answers depend, of
course, on $\mathcal{E}\rightarrow M$ we work with). If we choose the
alternative path II, the problem is to define the curvature of $\varepsilon (%
\mathcal{H}^{k})$ as we may not have any connection with its curvature at
our disposal. Of course, the solution is hinted by Propositions 4, 8 and
their proofs.

\begin{proposition}
The nonlinear PDE $\varepsilon (\mathcal{H}^{k})\subset J^{k+1}(\mathcal{E})$
has unique holonomic solutions with arbitrary inital conditions if and only
if the map $\varrho (\varepsilon (\mathcal{H}^{k}))\rightarrow \varepsilon (%
\mathcal{H}^{k})$ is surjective. In fact, there exists a vector bundle $%
F_{1} $ over $\varepsilon (\mathcal{H}^{k})$ and a section $\mathcal{R}$ of
this bundle satisfying the following diagram which is exact at the middle
term

\begin{equation}
\varrho (\varepsilon (\mathcal{H}^{k}))\overset{\pi _{k+1}^{k+2}}{%
\longrightarrow }\varepsilon (\mathcal{H}^{k})\overset{\mathcal{R}}{%
\longrightarrow }F_{1}
\end{equation}
\end{proposition}

So $\mathcal{R=}0$ if and only if $\varrho (\varepsilon (\mathcal{H}%
^{k}))\rightarrow \varepsilon (\mathcal{H}^{k})$ is onto.

\begin{definition}
$\mathcal{R}$ is the curvature of the phg $\varepsilon (\mathcal{H}%
^{k})\subset J^{k+1}(\mathcal{E}).$
\end{definition}

Proposition 10 is a special case of Theorem 1, pg. 95, [P2], where $%
\varepsilon (\mathcal{H}^{k})$ is replaced by a general nonlinear PDE $%
\mathcal{K}^{k+1}\subset J^{k+1}(\mathcal{E})$ and $\varrho (\varepsilon (%
\mathcal{H}^{k}))$ by its prolongation $\varrho (\mathcal{K}^{k+1}).$ The
proof of this general case is rather technical and not easy to verify even
on concrete examples, as the author states. Furthermore, $\mathcal{R}$ is
only the first of a sequence of curvatures $\mathcal{R}_{1},\mathcal{R}%
_{2},...,$ whose vanishing implies surjectivity of higher order
prolongations. These methods imply only formal integrability in smooth
category. Now a phg, as we remarked above, is a very special PDE and its
definition is motivated by the stabilization order of Klein geometries $%
(G,H) $ ([Or1], [Ol2]). For instance, the symbol of a phg vanishes, i.e., it
is of type zero and therefore trivially $k$-acyclic for any $k$ and
involutive. Therefore all the assumptions of Theorem 1 [P2] are satisfied by
a phg. Another simplification is that formal integrability can be replaced
by strong integrability as stated by Proposition 10 whose proof, as
expected, reduces to Frobenius Theorem.

It remains to check Proposition 10 on projective and conformal structures as
a preperation for higher order structures ([Or5]) .

\bigskip

\textbf{References}

\bigskip

[B] A.D.Blaom: Geometric structures as deformed infinitesimal symmetries,
Trams. Amer. Math. Soc., 358, 2006, 2651-71

[Gm] W.M.Goldman: Geometric Structures on Manifolds, AMS, Graduate Studies
in Math., 227, 2022

[Gs] H. Goldschmidt: Integrability criterion for systems of nonlinear
partial differential equations, J. Differential Geom. 1, (1969), 269-307

[KS] A.Kumpera, D.C.Spencer: Lie Equations, Ann. Math. Studies, 73,
Princeton Univ. Press, Princeton, New Jersey, 1972

[KLV] I.S.Krasil'shchik, V.V.Lychagin, A.M.Vinogradov: Geometry of Jet
Spaces and Nonlinear Partial Differential Equations, Gordon and Breach, 1986

[Ol1] P.J.Olver: Applications of Lie Groups to Differential Equations,
Springer, Berlin, 1986

[Ol2] \ \ \ \ \ \ " \ \ \ \ \ \ : Equivalence, Invariants and Symmetry,
Cambridge University Press, 1995

[Or1] E.H.Orta\c{c}gil: An Alternative Approach to Lie Groups and Geometric
Structures, OUP, 2018

[Or2] \ \ " \ \ \ : Curvature without connection, arXiv 2003.06593

[Or3] \ \ " \ \ : The mystery of torsion in differential geometry,
Researchgate

[Or4] \ \ " \ \ : Klein geometries of high order, arXiv 2211.02355

[Or5] \ \ " \ \ : Projective and conformal prehomogeneous geometries, in
progress

[P1] J.F.Pommaret : Systems of Partial Differential Equations and Lie
Pseudogroups, Gordon and Breach, London, New York, 1978

[P2] \ \ \ " \ \ \ : Partial Differential Equations and Group Theory, New
Perspectives for Applications, Kluwer Academic Publishers, Mathematics and
its Applications, Vol. 293, 1994

[S] D.C.Spencer: Overdetermined systems of partial\c{s} differential
equations, Bull. Amer. Math. Soc. 75 (1965), 1-114

[T] W.P.Thurston: Three-Dimensional Geometry and Topology, Vol.1, Princeton
University Press, 1997

\bigskip

ortacgile@gmail.com

\end{document}